\documentclass[11pt]{amsart}
\usepackage{amsmath,amssymb,amsfonts}
\usepackage[latin1]{inputenc}
\usepackage{epsfig,psfrag}
\numberwithin{equation}{section}
\newtheorem{thm}{Theorem}[section]
\newtheorem{aadef}[thm]{Definition}
\newtheorem{alem}[thm]{Lemma}
\newtheorem{aprop}[thm]{Proposition}
\newtheorem{acor}[thm]{Corollary}
\newtheorem{arem}[thm]{Remark}
\newtheorem{ax}[thm]{Example}\newenvironment{adem}[1][]%
   {\ \\ {\bf Proof #1. }}%
   {\hfill\mbox{\rule{2 true mm}{3 true mm}}}
   {\ \\{\bf Example #1. }}%
   {\hfill\mbox{\rule{2 true mm}{         3 true mm}}}
\newcommand{\R}{\mathbb{R}}
\renewcommand{\P}{\mathbb{P}}

\newcommand{\cT}{{\mathcal T}}
\newcommand{\E}{{\mathbb E}}

\newcommand{\N}{{\mathbb N}}

\newcommand{\lecx}{\le_{\textup{cx}}}
\newcommand{\mycomments}[1]{  }

\newcommand{\be}{ \begin{eqnarray*}  }
\newcommand{\ee}{ \end{eqnarray*}  }

\begin{document}
\title[Squared quadratic Wasserstein distance: optimal couplings and differentiability]{Squared quadratic Wasserstein distance : optimal couplings and Lions differentiability}
\author{Aur\'elien Alfonsi and Benjamin Jourdain}
\address{Universit\'e Paris-Est, Cermics (ENPC), INRIA, F-77455 Marne-la-Vall\'ee, France.}
\email{aurelien.alfonsi@enpc.fr, benjamin.jourdain@enpc.fr.}
\thanks{This research
 benefited from the support of the ``Chaire Risques Financiers'', Fondation du Risque.}
\date{\today}\maketitle
\begin{abstract}
  In this paper, we remark that any optimal coupling for the quadratic Wasserstein distance $W^2_2(\mu,\nu)$ between two probability measures $\mu$ and $\nu$ with finite second order moments on $\R^d$ is the composition of a martingale coupling with an optimal transport map ${\mathcal T}$. We check the existence of an optimal coupling in which this map gives the unique optimal coupling between $\mu$ and ${\mathcal T}\#\mu$. Next, we give a direct proof that $\sigma\mapsto W_2^2(\sigma,\nu)$ is differentiable at $\mu$ in the Lions~\cite{Lions} sense iff there is a unique optimal coupling between $\mu$ and $\nu$ and this coupling is given by a map. It was known combining results by Ambrosio, Gigli and Savar\'e~\cite{AGS} and Ambrosio and Gangbo~\cite{AmGa} that, under the latter condition, geometric differentiability holds. Moreover, the two notions of differentiability are equivalent according to the recent paper~\cite{GaTu} of Gangbo and Tudorascu. 
  Besides, we give a self-contained probabilistic proof that mere Fr\'echet differentiability of a law invariant function $F$ on $L^2(\Omega,\P;\R^d)$ is enough for the Fr\'echet differential at $X$ to be a measurable function of $X$. 
\end{abstract}

 \section*{Introduction}
 In this paper, we are interested in the structure of optimal couplings for the squared quadratic Wasserstein distance $W^2_2(\mu,\nu)$ between $\mu$ and $\nu$ in the set ${\mathcal P}_2(\R^d)$ of probability measures with finite second order moments on $\R^d$, and in the differentiability of $W^2_2(\mu,\nu)$ with respect to $\mu$. By definition, $W^2_2(\mu,\nu)=\inf_{\pi\in\Pi(\mu,\nu)} \int|y-x|^2\pi(dx,dy)$ where $\Pi(\mu,\nu)$ denotes the set of coupling measures on $\R^d\times\R^d$ with first and second marginals respectively equal to $\mu$ and $\nu$ and $|.|$ denotes the Euclidean norm on $\R^d$. There always exists an optimal coupling and we denote by $\Pi^{opt}(\mu,\nu)$ the set of optimal couplings. 
 According to~\cite{Gigli}, there exists only one $W_2$-optimal coupling $\pi$ between $\mu$ and each $\nu\in{\mathcal P}_2(\R^d)$ and this coupling is given by a map $T$ (i.e. $\pi=(I_d,T)\#\mu$ where $I_d$ denotes the identity function on $\R^d$) iff $\mu$ gives $0$ mass to the $c-c$ hypersurfaces of dimension $d-1$. Even when $\mu$ does not satisfy this condition which is implied by absolute continuity with respect to the Lebesgue measure, according to Proposition 5.13~\cite{CaDe}, if $\varphi:\R^d\to\R$ is a $C^2$ strictly convex function such that $\int_{\R^d}|\nabla\varphi(x)|^2\mu(dx)<\infty$, then there is a unique $W_2$-optimal coupling between $\mu$ and $\nu=\nabla\varphi\#\mu$ and this coupling is given by the map $\nabla\varphi$. But there also exist measures $\nu\in{\mathcal P}_2(\R^d)$ such that either the unique optimal coupling (uniqueness holds in dimension $d=1$ for instance) is not given by a map or there exist distinct optimal couplings. In the latter case, any strictly convex combination of these couplings is an optimal coupling which is not given by a map.

 In Section~\ref{secstruccoupl}, we study optimal couplings $\pi$ which are not given by a map. By disintegration, $\pi(dx,dy)=\mu(dx)k(x,dy)$ for some Markov kernel $k$ on $\R^d$ (which is $\mu(dx)$ a.e. unique). Setting ${\mathcal T}(x)=\int_{\R^d}yk(x,dy)$ and using the bias-variance decomposition under the kernel $k$, we obtain that $\pi$ is the composition of a martingale coupling between ${\mathcal T}\#\mu$ and $\nu$ with the map ${\mathcal T}$ which gives a $W_2$-optimal coupling between $\mu$ and ${\mathcal T}\#\mu$.  Note that couplings of this form have recently been studied by Gozlan and Juillet~\cite{GoJu} when considering the barycentric optimal cost problem.  For $\phi:\R^d\to\R$ a strictly convex function such that $\int_{\R^d}\phi(y)\nu(dy)<\infty$, by minimizing $\int_{\R^d}\phi({\mathcal T}(x))\mu(dx)$ over the $W_2$-optimal couplings between $\mu$ and $\nu$, we obtain optimal couplings such that the associated map ${\mathcal T}_\phi$ gives the only optimal coupling between $\mu$ and ${\mathcal T}_\phi\#\mu$. There is a unique such coupling when $\phi(x)=|x|^2$.
 
 In Section~\ref{secldif}, we are interested in the differentiability of $W^2_2(\mu,\nu)$ in the Lions sense with respect to $\mu$.  Gangbo and Tudorascu have recently proved in Corollary 3.22 \cite{GaTu} that the Lions differentiability~\cite{Lions} of a function $f:\mathcal{P}_2(\R^d)\rightarrow \R$  is equivalent to the geometric differentiability and that the Fr\'echet derivative of the lift at $X\sim\mu$ is then given by $\nabla_\mu f(X)$ where $\nabla_\mu f\in L^2(\R^d,\mu;\R^d)$ is the geometric (or Wasserstein) gradient of $f$ at $\mu$. While the lifted space that they consider is the ball centered at the origin of unit volume in~$\R^d$ endowed with the Lebesgue measure, the result can be transferred to any atomless lifted space by considering an almost isomorphism between those spaces\footnote{We thank one of the referees for pointing out this argument to us.}. Theorem 10.2.6~\cite{AGS} states that $\sigma\mapsto W_2^2(\sigma,\nu)$ is subdifferentiable in the geometric sense at $\mu$ when $\Pi^{opt}(\mu,\nu)=\{(I_d,T)\#\mu\}$ for some measurable transport map $T:\R^d\to\R^d$. On the other hand, Proposition 4.3~\cite{AmGa} states that $\sigma\mapsto W_2^2(\sigma,\nu)$ is always superdifferentiable in the geometric sense at $\mu$ with $x\mapsto 2\left(x-\int_{\R^d}yk(x,dy)\right)$ belonging to the superdifferential for each Markov kernel $k$ on $\R^d$ such that $\mu(dx)k(x,dy)\in\Pi^{opt}(\mu,\nu)$. Since geometric differentiability amounts to simultaneous geometric sub- and superdifferentiability, as soon as $\Pi^{opt}(\mu,\nu)=\{(I_d,T)\#\mu\}$, then $\sigma\mapsto W_2^2(\sigma,\nu)$ is differentiable in the geometric sense at $\mu$. On the other hand, geometric differentiability implies that the geometric sub- and superdifferential considered as subsets of $L^2(\R^d,\mu;\R^d)$ coincide and contain one element only (see for instance Proposition 5.63 \cite{CaDe}). The fact that the quotient of $\{x\mapsto\int_{\R^d}yk(x,dy):\mu(dx)k(x,dy)\in\Pi^{opt}(\mu,\nu)\}$ for the $\mu(dx)$ a.e. equality is a singleton is therefore necessary for the geometric differentiability of $\sigma\mapsto W_2^2(\sigma,\nu)$ to hold at $\mu$.\\
We prove that $\sigma\mapsto W_2^2(\sigma,\nu)$ is differentiable at $\mu$ in the Lions sense iff $\Pi^{opt}(\mu,\nu)=\{(I_d,T)\#\mu\}$. We give a direct probabilistic proof of the sufficient condition which also follows from the just mentionned results. To prove the necessary condition, we use that the Fr\'echet differentiability at $X\sim\mu$ of the lift on an atomless probability space is enough for the Fr\'echet derivative at $X$ to be a.s. equal to a measurable function of $X$, a consequence of~\cite{GaTu} that we show again using simple probabilistic arguments. Let us emphasize that the quotient of $\{x\mapsto\int_{\R^d}yk(x,dy):\mu(dx)k(x,dy)\in\Pi^{opt}(\mu,\nu)\}$ for the $\mu(dx)$ a.e. equality may be a singleton while $\Pi^{opt}(\mu,\nu)$ is not equal to $\{(I_d,T)\#\mu\}$ for some measurable map $T:\R^d\to\R^d$ (see, in dimension $d=1$, Remark \ref{remexnondif} below).

{\bf Acknowledgements :} We thank Pierre Cardaliaguet for pointing out the reference \cite{GaTu} to us.
\section{Structure of quadratic Wasserstein optimal couplings}\label{secstruccoupl}

In this section, we are interested in characterizing the set \begin{align*}
  \Pi^{opt}(\mu,\nu)=\{ \pi(dx,dy)\in {\mathcal P}_2(\R^d\times \R^d) : &\mu(dx)=\int_{y\in\R^d} \pi(dx,dy), \nu(dy)=\int_{x\in\R^d} \pi(dx,dy) \\&\text{ and }W_2^2(\mu,\nu)=\int_{\R^d \times \R^d}|y-x|^2 \pi(dx,dy) \}.
                                                             \end{align*} of optimal couplings between two probability measures $\mu,\nu \in {\mathcal P}_2(\R^d)$ for the quadratic cost. This set is not empty : see e.g.~\cite{AGS} p.~133.
                                                             
The refined version of the Brenier theorem in~\cite{Gigli} ensures that $\Pi^{opt}(\mu,\nu)$ contains a single element $(I_d,T)\#\mu$ which is given by a measurable transport map $T:\R^d\rightarrow \R^d$ for each $\nu\in{\mathcal P}_2(\R^d)$ iff $\mu$ does not give mass to the $c-c$ hypersurfaces parametrized by an index $i\in\{0,\hdots,d-1\}$ and two convex functions $f$ and $g$ from $\R^{d-1}$ to $\R$ : $$\{(x_1,\hdots,x_{i},f(x)-g(x),x_{i+1},\hdots,x_{d-1}):x=(x_1,\hdots,x_{d-1})\in\R^{d-1}\}.$$ 
The next lemma deals with the case where $\Pi^{opt}(\mu,\nu)\not=\{(I_d,T)\#\mu\}$ for some measurable transport map. 

\begin{alem}\label{cor_utile}Let $\mu,\nu \in \mathcal{P}_2(\R^d)$. One of the two conditions holds:
  \begin{itemize}
    \item $\Pi^{opt}(\mu,\nu)=\{(I_d,T)\#\mu\}$ for some measurable transport map $T:\R^d\rightarrow \R^d$, \item $\exists \mu(dx)k(x,dy)\in \Pi^{opt}(\mu,\nu)$ such that $\int_{\R^d\times\R^d}|y-\int_{\R^d}zk(x,dz)|^2k(x,dy)\mu(dx)>0$.
    \end{itemize}
    Moreover, if any coupling in $\Pi^{opt}(\mu,\nu)$ is given by a map i.e. writes $(I_d,T)\#\mu$ for some measurable function $T:\R^d\to\R^d$, then $\Pi^{opt}(\mu,\nu)$ is a singleton.
\end{alem}
\begin{adem}
If the set $\Pi^{opt}(\mu,\nu)$ has a single element $\mu(dx)k(x,dy)$, defining $\mathcal{T}(x)=\int_{\R^d}yk(x,dy)$ we either have  $\int_{\R^d\times\R^d}|y-\mathcal{T}(x)|^2k(x,dy)\mu(dx)>0$ or $\mu(dx)k(x,dy)=\mu(dx)\delta_{\mathcal{T}(x)}(dy)$. Otherwise, we can pick two distinct elements $k_1,k_2\in\Pi^{opt}(\mu,\nu)$ and $k(x,dy)=\frac{1}{2}\left(k_1(x,dy)+k_2(x,dy)\right)$ is such that $\mu(dx)k(x,dy)\in\Pi^{opt}(\mu,\nu)$ and $\int_{\R^d\times\R^d}|y-\int_{\R^d}zk(x,dz)|^2k(x,dy)\mu(dx)>0$. The second statement easily follows.
\end{adem}

Remarking that if $\nu$ is the Dirac mass at $x\in\R^d$ and $\nu_\varepsilon$ the uniform distribution on the ball centered at $x$ with radius $\varepsilon$, then $W_2(\nu,\nu_\varepsilon)\le \varepsilon$, we deduce from the next proposition that for any $\mu,\nu \in {\mathcal P}_2(\R^d)$, we can always find $\mu_\varepsilon,\nu_\varepsilon \in {\mathcal P}_2(\R^d)$ such that $W_2(\mu,\mu_\varepsilon)\le \varepsilon$, $W_2(\nu,\nu_\varepsilon)\le \varepsilon$ and 
$\exists \mu_\varepsilon(dx)k_\varepsilon(x,dy)\in \Pi^{opt}(\mu_\varepsilon,\nu_\varepsilon)$ such that $\int_{\R^d\times\R^d}|y-\int_{\R^d}zk_\varepsilon(x,dz)|^2k_\varepsilon(x,dy)\mu_\varepsilon(dx)>0$. 
\begin{aprop}\label{propnondifvois}
   Assume that $\nu\in{\mathcal P}_2(\R^d)$ is not a Dirac mass. Then for all $\mu\in{\mathcal P}_2(\R^d)$, there exists a sequence $(\mu_n)_{n}$ of elements of ${\mathcal P}_2(\R^d)$ such that $\lim_{n\to\infty}W_2(\mu_n,\mu)=0$ and for each $n$, there does not exist $T_n:\R^d\rightarrow \R^d$ measurable such that $\Pi^{opt}(\mu_n,\nu)=\{(I_d,T_n)\#\mu_n\}$.
\end{aprop}
\begin{adem}
  Let $(X_i)_{i\ge 1}$ be an i.i.d. sequence of random variables with law $\mu$, and $(Y_i)_{i\ge 1}$ an independent  i.i.d. sequence of uniform random variables on the unit ball $\{x\in \R^d, |x|\le 1\}$. We set $\tilde{\mu}_n=\frac{1}{n}\sum_{i=1}^n\delta_{X_i}$ the empirical measure and $\mu_n=\frac{1}{n}\sum_{i=1}^n\delta_{X_i+Y_i/n}$. By construction, we have $W_2^2(\mu_n,\tilde{\mu}_n)\le \frac{1}{n} \sum_{i=1}^n |Y_i/n|^2 \le 1/n^2$ and $\P(\exists i \not = j, X_i+Y_i/n=X_j+Y_j/n)=0$, which means that a.s. for each $n\in \N^*$, $\mu_n$ weights a.s. exactly $n$ points.
  The law of large numbers gives the almost sure weak convergence of $  \tilde{\mu}_n$ towards $\mu$ and the almost sure convergence of $\frac 1n \sum_{i=1}^n |X_i|^2$ to $\E[|X_1|^2]$. Proposition 7.1.5 in~\cite{AGS} ensures that $W_2(\tilde{\mu}_n,\mu)\underset{n \rightarrow +\infty}\rightarrow 0$ almost surely. By the triangle inequality, we get $W_2(\mu_n,\mu)\underset{n \rightarrow +\infty}\rightarrow 0$ almost surely.

  Now, we consider $(p_n)_{n\ge 1}$ the increasing sequence of prime numbers.  Suppose that $\exists n_0 \in \N^*$, such that $T\# \mu_{p_{n_0}}=\nu$. Then, $\nu$  weights at most $p_{n_0}$ points and the masses are equal to $k/p_{n_0}$ with $1\le k\le p_{n_0}-1$ since $\nu$ is not a Dirac mass. Then, if we had $T\# \mu_{p_{n}}=\nu$ for some $n>n_0$, we would have $k/p_{n_0}=k'/p_n$ with  $1\le k'\le p_{n}-1$. This would imply that $p_{n_0}$ divides $kp_n$ and thus $k$, which is impossible since $1\le k\le p_{n_0}-1$. Thus, there is at most one $n_0 \in \N^*$ such that  there is a transport map $T_{n_0}$ satisfying $T_{n_0}\# \mu_{p_{n_0}}=\nu$.
\end{adem}

Let us now give a necessary and sufficient condition for the existence of an optimal transport map in dimension~$d=1$. We denote $F_\eta(x)=\eta((-\infty,x])$ and $F_\eta^{-1}(u)=\inf\{x\in\R:\eta((-\infty,x])\ge u\}$ the cumulative distribution function and the quantile function of a probability measure $\eta$ on $\R$. For $\mu,\nu\in{\mathcal P}_2(\R)$, by Theorem 2.9 in~\cite{Santambrogio}, the only element of $\Pi^{opt}(\mu,\nu)$ is the image of the Lebesgue measure on $[0,1]$ by $(F_\mu^{-1},F_\nu^{-1})$. The next lemma characterizes the case when this coupling is given by a map.
\begin{alem}\label{lem_map1d}
   Let $\mu,\nu\in{\mathcal P}_2(\R)$. There exists $T\in L^2(\R,\mu;\R)$  such that $\Pi^{opt}(\mu,\nu)=\{(I_1,T)\# \mu\}$ iff for all $x\in \R$ such that $\mu(\{x\})>0$, $F_\nu^{-1}$ is constant on $(F_\mu(x-),F_\mu(x)]$. Then, the unique optimal transport map is $T(x)=F_\nu^{-1}(F_\mu(x))$.
\end{alem}
\begin{arem}\label{remexnondif}
   When $F_\nu^{-1}$ is not constant on $(F_\mu(x-),F_\mu(x)]$ for some $x\in\R$ such that $\mu(\{x\})>0$, then $\Pi^{opt}(\mu,\nu)$ is not equal to $\{(I_1,T)\# \mu\}$ for some measurable map $T:\R\to\R$ while, since $\Pi^{opt}(\mu,\nu)$ is a singleton, the quotient of $\{x\mapsto\int_{\R^d}yk(x,dy):\mu(dx)k(x,dy)\in\Pi^{opt}(\mu,\nu)\}$ for the $\mu(dx)$ a.e. equality is a singleton.
\end{arem}

\begin{adem}
 Let $X\sim \mu$ and $U$ be an independent random variable uniform on $[0,1]$. The random variable $V=F_\mu(X-)+U(F_\mu(X)-F_\mu(X-))$ is such that $\P(\{F_\mu(X-)<V\le F_\mu(X)\} \cup \{F_\mu(X-)=V= F_\mu(X)\} )=1$.  This is an uniform random variable on $[0,1]$: for $u\in(0,1)$, $u\in [F_\mu(x-),F_\mu(x)]$ for some $x\in \R$ and $\P(V\le u)=\P(X<x)+\P\left(X=x,U\le \frac{u-F_\mu(x-)}{F_\mu(x)-F_\mu(x-)}\right)=u$ since $X$ is independent of $U$. Since $F_\mu^{-1}(V)= X$ for $V\in (F_\mu(X-),F_\mu(X)]$ and  $F_\mu^{-1}(V)\le X$ for $V=F_\mu(X-)=F_\mu(X)$, we have $F_\mu^{-1}(V)\le X$ a.s.. Since $F_\mu^{-1}(V)$ and $X$  have the same law, we necessarily have $F_\mu^{-1}(V)= X$ a.s.. By the inverse transform sampling, $F^{-1}_\nu(V)$ is distributed according to $\nu$. Let us assume that  $F_\nu^{-1}$ is constant on $(F_\mu(x-),F_\mu(x)]$ for all $x\in \R$ such that $\mu(\{x\})>0$. Then $F^{-1}_\nu(V)= F^{-1}_\nu(F_\mu(X))$ a.s.,  $F_\nu^{-1}\circ F_\mu\#\mu=\nu$ and
          $$\int_0^1 (F^{-1}_\mu(v)-F^{-1}_\nu(v))^2dv = \E[(X- F^{-1}_\nu(F_\mu(X)))^2]=\int_{\R}(x- F^{-1}_\nu(F_\mu(x)))^2\mu(dx).$$
 Hence $T(x)=F_\nu^{-1}(F_\mu(x))$ is an optimal transport map. 
Conversely, if $T$ is an optimal transport map such that $T\#\mu=\nu$, we have $T(F^{-1}_\mu(v))=F^{-1}_\nu(v)$, $dv$-a.e. For $x\in \R$ such that $\mu(\{x\})>0$, $F_\mu^{-1}$ is constant on $(F_\mu(x-),F_\mu(x)]$, and therefore $F_\nu^{-1}$ is necessarily constant on $(F_\mu(x-),F_\mu(x)]$. 
\end{adem}

\begin{arem} Lemma~\ref{lem_map1d} still holds true  for $\mu,\nu$ probability measures on $\R$ with finite moments of order $\rho\ge 1$, and a transport cost $c(x,y)=h(|y-x|)$, with $h:\R_+\rightarrow \R$ strictly convex such that $\exists C<\infty,\;\forall x\in\R,\;h(|x|)\le C(1+|x|^\rho)$. The same proof applies since, by Theorem 2.9 in~\cite{Santambrogio}, the only optimal coupling for such a cost is the image of the Lebesgue measure on $[0,1]$ by $(F_\mu^{-1},F_\nu^{-1})$.
\end{arem}

The next proposition, which is one of the main results of this section, shows that any $W_2$-optimal coupling can be written as the composition of a transport map and a martingale kernel i.e. a Markov kernel $k$ such that for all $x\in\R^d$, $\int_{\R^d}|y|k(x,dy)<\infty$ and $\int_{\R^d}yk(x,dy)=x$. Let us now give the definition of the convex order on probability measures before recalling its link with the existence of martingale couplings.
\begin{aadef}
  Let $\eta,\nu$ be two probability measures on $\R^d$. We say that $\eta$ is smaller than $\nu$ in the convex order and write $\eta\lecx \nu$ if for each convex function $\phi:\R^d\to\R$ such that the integrals make sense,
  $$\int_{\R^d}\phi(x)\eta(dx)\le\int_{\R^d}\phi(y)\nu(dy).$$
\end{aadef}
Notice that since a convex function $\phi$ on $\R^d$ is bounded from below by an affine function, for a probability measure $\eta$ on $\R^d$ with finite first order moment (and in particular for $\eta\in{\mathcal P}_2(\R^d)$), $\int_{\R^d}\phi(x)\eta(dx)$ always makes sense possibly equal to $+\infty$.

Theorem 8 in Strassen~\cite{Strassen} ensures that, when $\int_{\R^d}|y|\nu(dy)<\infty$, $\eta\lecx \nu$ iff there exists a martingale Markov kernel $k$ such that $\eta(dx)k(x,dy)\in\Pi(\eta,\nu)$.

\begin{aprop}\label{lemstructoplan}Let $\mu,\nu \in{\mathcal P}_2(\R^d)$, $\mu(dx)k(x,dy)\in \Pi^{opt}(\mu,\nu)$, ${\mathcal T}(x)=\int_{\R^d}yk(x,dy)$ and $\eta={\mathcal T}\#\mu$. Then $\eta\lecx\nu$, \begin{equation}
   W_2^2(\mu,\nu)=W_2^2(\mu,\eta)+\int_{\R^d}|y|^2\nu(dy)-\int_{\R^d}|z|^2\eta(dz)\label{decompw2}
 \end{equation}
  and $(I_d,\mathcal{T})\#\mu\in\Pi^{opt}(\mu,\eta)$.
  
 On the other hand, if $\eta\lecx\nu$ is such that \eqref{decompw2} holds, then combining $\mu(dx)q(x,dz)\in\Pi^{opt}(\mu,\eta)$ with any martingale coupling $\eta(dz)m(z,dy)$ between $\eta$ and $\nu$, we obtain a $W_2$-optimal coupling $\mu(dx)qm(x,dy)$ (where, as usual, $qm(x,dy)=\int_{z\in\R^d}q(x,dz)m(z,dy)$) between $\mu$ and $\nu$.
\end{aprop}
The first part of this proposition is also a consequence of Theorem 12.4.4 in~\cite{AGS} : the barycentric projection of $\mu(x)k(x,dy)$ is precisely $(I_d,\mathcal{T})\# \mu$. Here, we present this result with a probabilistic fashion. For $\mu(dx)k(x,dy)$ as in the first statement and $(X,Y)\sim\mu(dx)k(x,dy)$, by definition of ${\mathcal T}$, $\E[Y|X]={\mathcal T}(X)$ a.s. so that $\E[Y|{\mathcal T}(X)]={\mathcal T}(X)$ a.s. and this optimal coupling is the composition of the martingale coupling given by the law of $({\mathcal T}(X),Y)$ and the transport map ${\mathcal T}$. Notice that since it relies on the bias-variance decomposition, this structure of optimal couplings does not seem to generalize to other Wasserstein distances $W_\rho(\mu,\nu)=\left( \inf_{\pi\in\Pi(\mu,\nu)} \int|y-x|^\rho\pi(dx,dy) \right)^{1/\rho}$, $\rho \in [1,\infty) \setminus \{2\}$. Nevertheless, Gozlan and Juillet~\cite{GoJu} have recently obtained optimal couplings that are the composition of a martingale coupling  and a deterministic transport map by considering the barycentric optimal cost problem, which consists in minimizing for a given cost function $\theta:\R^d \rightarrow \R_+$ the quantity $\int_{\R^d} \theta(x-\int_{\R^d}yk(x,dy))\mu(dx)$ among all couplings $\mu(dx)k(x,dy)$ between $\mu$ and $\nu$.

\begin{adem}
  Let us first prove the second statement. Let $\eta\lecx\nu$, $q$ be a Markov kernel such that $\mu(dx)q(x,dz)\in\Pi^{opt}(\mu,\eta)$ and $m$ be any martingale kernel such that $\eta m=\nu$. Then $\mu(dx)qm(x,dy)$ is a coupling between $\mu$ and $\nu$ such that
  \begin{align}
   W_2^2(\mu,\nu)&\le \int_{\R^d\times\R^d}|y-x|^2\mu(dx)qm(x,dy)=\int_{\R^d\times\R^d\times\R^d}|y-z+z-x|^2\mu(dx)q(x,dz)m(z,dy)\notag\\&=\int_{\R^d\times\R^d}|y-z|^2\eta(dz)m(z,dy)+\int_{\R^d\times\R^d}|z-x|^2\mu(dx)q(x,dz)\notag\\&=\int_{\R^d}|y|^2\nu(dy)-\int_{\R^d}|z|^2\eta(dz)+W_2^2(\mu,\eta)\label{decompo}
  \end{align}
  where we used the variance-bias decomposition under the martingale kernel $m$ for the third equality. Hence, if \eqref{decompw2} holds, then $\mu(dx)qm(x,dy)\in\Pi^{opt}(\mu,\nu)$.

  Let now $\mu(dx)k(x,dy)\in\Pi^{opt}(\mu,\nu)$, ${\mathcal T}(x)=\int_{\R^d}yk(x,dy)$ and $\eta={\mathcal T}\#\mu$. Jensen's inequality immediately gives $\eta\lecx \nu$ and thus $\eta \in {\mathcal P}_2(\R^d)$. We have
  \begin{align*}
    W_2^2(\mu,\nu)&=\int_{\R^d}\int_{\R^d}|y-\mathcal{T}(x)+\mathcal{T}(x)-x|^2\mu(dx)k(x,dy)\\
    &=\int_{\R^d}\int_{\R^d}|y-\mathcal{T}(x)|^2\mu(dx)k(x,dy) + \int_{\R^d}|\mathcal{T}(x)-x|^2\mu(dx)\\
    &=\int_{\R^d}\int_{\R^d}(|y|^2-|\mathcal{T}(x)|^2)\mu(dx)k(x,dy) + \int_{\R^d}|\mathcal{T}(x)-x|^2\mu(dx)\\
    &=\int_{\R^d}|y|^2\nu(dy)-\int_{\R^d}|z|^2\eta(dz) + \int_{\R^d}|\mathcal{T}(x)-x|^2\mu(dx) ,
  \end{align*}
  where we used the variance-bias decomposition with respect to $k(x,.)$ for the second equality. With \eqref{decompo}, we deduce that $\int_{\R^d}|\mathcal{T}(x)-x|^2\mu(dx)\le W_2^2(\mu,\eta)$ and $\mathcal{T}$ is a $W_2$-optimal  transport map between $\mu$ and $\eta$. 
\end{adem}
  
 For $\mu,\nu\in {\mathcal P}_2(\R^d)$,  let us define the sets
    \begin{align*}{\mathcal I}_\mu^\nu&=\{\eta\in{\mathcal P}_2(\R^d):\eta\lecx\nu\mbox{ and }W_2^2(\mu,\nu)=W_2^2(\mu,\eta)+\int_{\R^d}|y|^2\nu(dy)-\int_{\R^d}|z|^2\eta(dz)\}, \\
 \tilde{{\mathcal I}}_\mu^\nu&=\left\{{\mathcal T}\#\mu : \exists \mu(dx)k(x,dy)\in \Pi^{opt}(\mu,\nu),  {\mathcal T}(x)=\int_{\R^d}yk(x,dy) \right\} .  \end{align*}By Proposition~\ref{lemstructoplan}, we have $\tilde{{\mathcal I}}_\mu^\nu \subset {\mathcal I}_\mu^\nu$ and $ \tilde{{\mathcal I}}_\mu^\nu \not = \emptyset $ since $\Pi^{opt}(\mu,\nu) \neq \emptyset$. Moreover, there exists an optimal transport map between $\mu$ and any element of $\tilde{{\mathcal I}}_\mu^\nu$. The measure ${\mathcal T}\#\mu$ associated with an optimal coupling in $\Pi^{opt}(\mu,\nu)$ is possibly equal to $\nu$, which always belongs to ${\mathcal I}_\mu^\nu$. \begin{alem}\label{cor_Imunu}
Let $\mu,\nu \in \mathcal{P}_2(\R^d)$. If $\eta \in \mathcal{I}_\mu^\nu$, then for any $\tilde{\eta}$ such that $\eta\lecx \tilde{\eta} \lecx \nu$, $\tilde{\eta}\in   \mathcal{I}_\mu^\nu$ and $\eta \in \mathcal{I}_\mu^{\tilde{\eta}}$. Moreover, $\mathcal{I}_\mu^\nu=\{\eta\in{\mathcal P}_2(\R^d):\exists\tilde\eta \in\tilde{\mathcal{I}}_\mu^\nu,\tilde\eta\lecx\eta\lecx\nu\}$. Last, the set $\mathcal{I}_\mu^\nu$ is convex. 
\end{alem}
\begin{adem}
  Let $\eta \in \mathcal{I}_\mu^\nu$ and $\tilde{\eta}$ be such that $\eta\lecx \tilde{\eta} \lecx \nu$. We have \begin{equation}\label{decompo_wass}
    W_2^2(\mu,\nu)=W_2^2(\mu,\eta)+\int_{\R^d}|y|^2\nu(dy)-\int_{\R^d}|\tilde{z}|^2\tilde{\eta}(d\tilde{z})+\int_{\R^d}|\tilde{z}|^2\tilde{\eta}(d\tilde{z})-\int_{\R^d}|z|^2\eta(dz).
  \end{equation}
  Now, we consider $\mu(dx)k(x,dz)\in \Pi^{opt}(\mu,\eta)$ and $\eta(dz)m(z,d\tilde{z})$ a martingale coupling between $\eta$ and $\tilde{\eta}$. Then,  $W_2^2(\mu,\tilde{\eta})\le \int_{(\R^d)^3}|\tilde{z}-z+z-x|^2\mu(dx)k(x,dz)m(z,d\tilde{z})=W_2^2(\mu,\eta)+\int_{\R^d}|\tilde{z}|^2\tilde{\eta}(d\tilde{z})-\int_{\R^d}|z|^2\eta(dz)$. This inequality cannot be strict: otherwise, by combining an optimal coupling between $\mu$ and $\tilde{\eta}$ and a martingale coupling between $\tilde{\eta}$ and $\nu$, we would contradict~\eqref{decompo_wass}. The equality gives $\eta \in \mathcal{I}_\mu^{\tilde{\eta}}$ and $\tilde{\eta}\in   \mathcal{I}_\mu^\nu$ by using~\eqref{decompo_wass}.

  If $\tilde\eta \in\tilde{\mathcal{I}}_\mu^\nu$, since $\tilde{\mathcal{I}}_\mu^\nu\subset{\mathcal{I}}_\mu^\nu$, by the first statement, each probability measure $\eta$ such that $\tilde\eta\lecx\eta\lecx\nu$ belongs to ${\mathcal{I}}_\mu^\nu$. Hence $\{\eta\in{\mathcal P}_2(\R^d):\exists\tilde\eta \in\tilde{\mathcal{I}}_\mu^\nu,\tilde\eta\lecx\eta\lecx\nu\}\subset\mathcal{I}_\mu^\nu$. On the other hand, for $\eta\in{\mathcal{I}}_\mu^\nu$, $\mu(dx)q(x,dz)\in\Pi^{opt}(\mu,\eta)$ and a martingale coupling $\eta(dz)m(z,dy)$ between $\eta$ and $\nu$, we have $\mu(dx)qm(x,dy)\in\Pi^{opt}(\mu,\nu)$, by the second assertion in Proposition~\ref{lemstructoplan}. Since, by the martingale property, $\int_{\R^d}yqm(x,dy)=\int_{\R^d}\int_{\R^d}ym(z,dy)q(x,dz)=\int_{\R^d}zq(x,dz)$ setting ${\mathcal T}(x)=\int_{\R^d}zq(x,dz)$, we have ${\mathcal T}\#\mu\in\tilde{\mathcal{I}}_\mu^\nu$, by the first assertion in Proposition~\ref{lemstructoplan}. Since ${\mathcal T}\#\mu\lecx\eta$, we conclude that $\mathcal{I}_\mu^\nu\subset\{\eta\in{\mathcal P}_2(\R^d):\exists\tilde\eta \in\tilde{\mathcal{I}}_\mu^\nu,\tilde\eta\lecx\eta\lecx\nu\}$.

 Last, let us consider $\eta_1,\eta_2 \in \mathcal{I}_\mu^\nu$ and $\lambda \in (0,1)$. Using a convex combination of couplings in $\Pi^{opt}(\mu,\eta_1)$ and $\Pi^{opt}(\mu,\eta_2)$, we obtain that $W^2_2(\mu,\lambda \eta_1+(1-\lambda)\eta_2)\le \lambda W_2^2(\mu,\eta_1)+(1-\lambda)W_2^2(\mu,\eta_2)$. Since $\eta_1,\eta_2 \in \mathcal{I}_\mu^\nu$, we deduce that
 $$W_2^2(\mu,\nu)\ge W_2^2(\mu,\lambda \eta_1+(1-\lambda)\eta_2)+\int_{\R^d}|y|^2\nu(dy)-\int_{\R^d}|z|^2(\lambda \eta_1+(1-\lambda)\eta_2)(dz).$$
 Since $\lambda \eta_1+(1-\lambda)\eta_2 \lecx \nu$, there exists a martingale coupling between $\lambda \eta_1+(1-\lambda)\eta_2$ and $\nu$. Composing it with an element of $\Pi^{opt}(\mu,\lambda \eta_1+(1-\lambda)\eta_2)$, we obtain a coupling between $\mu$ and $\nu$ which ensures that
$$W_2^2(\mu,\nu)\le W_2^2(\mu,\lambda \eta_1+(1-\lambda)\eta_2)+\int_{\R^d}|y|^2\nu(dy)-\int_{\R^d}|z|^2(\lambda \eta_1+(1-\lambda)\eta_2)(dz).$$ Hence $\lambda \eta_1+(1-\lambda)\eta_2 \in \mathcal{I}_\mu^\nu$.
\end{adem}

In dimension $d=1$, since $\Pi^{opt}(\mu,\nu)$ is a singleton, we can specify the sets $\mathcal{I}_\mu^\nu $ and $\tilde{\mathcal{I}}_\mu^\nu$.

\begin{aprop} \label{propopt1d}Let $\mu,\nu\in{\mathcal P}_2(\R)$ and \begin{equation}
   \cT(x)=\int_0^1F_{\nu}^{-1}(F_{\mu}(x-)+u[F_{\mu}(x)-F_{\mu}(x-)] )du\label{defctd1}.
 \end{equation}
 We have $\tilde{\mathcal{I}}_\mu^\nu =\{{\mathcal T} \# \mu\}$ and ${\mathcal{I}}_\mu^\nu =\{\eta\in{\mathcal P}_2(\R):{\mathcal T} \# \mu\lecx\eta\lecx\nu\}$.  Moreover, $\Pi^{opt}(\mu,{\mathcal T}\#\mu)=\{(I_1,{\mathcal T})\#\mu\}$ and there is a unique martingale coupling between ${\mathcal T}\#\mu$ and $\nu$ and it is $W_2$-optimal.
\end{aprop}
\begin{adem}
By the second assertion in Lemma~\ref{cor_Imunu}, the characterization of ${\mathcal{I}}_\mu^\nu$ easily follows from the one of $\tilde{\mathcal{I}}_\mu^\nu$, which, with the definition of $\tilde{\mathcal{I}}_\mu^\nu$, the first statement in Proposition~\ref{lemstructoplan} and the uniqueness of the optimal coupling in dimension $d=1$, also implies that $\Pi^{opt}(\mu,{\mathcal T}\#\mu)=\{(I_1,{\mathcal T})\#\mu\}$. Let $U, U'$ be two independent uniform random variables on $[0,1]$. We define
   \begin{equation}\label{defV}
     V=F_{\mu}(F_{\mu}^{-1}(U)-)+U'[F_{\mu}(F_{\mu}^{-1}(U))-F_{\mu}(F_{\mu}^{-1}(U)-)],
   \end{equation}
and have by construction
\begin{equation}\label{property_V}
  F_{\mu}^{-1}(V)=F_{\mu}^{-1}(U)\mbox{ a.s.}.
\end{equation}

For $u\in(0,1)$, $u\in [F_\mu(x-),F_\mu(x)]$ for some $x\in \R$ and $$\P(V\le u)=\P(F_\mu^{-1}(U)<x)+\P\left(F_\mu^{-1}(U)=x,U'\le \frac{u-F_\mu(x-)}{F_\mu(x)-F_\mu(x-)}\right)=u$$ since $U'$ is independent of $U$.
Hence $V$ is uniformly distributed on $[0,1]$.   According to Theorem 2.9~\cite{Santambrogio}, the law of $(F_\mu^{-1}(V),F_\nu^{-1}(V))$ is the unique element of $\Pi^{opt}(\mu,\nu)$. From~\eqref{defV}, we get $\E[F_{\nu}^{-1}(V)|U]=\cT(F_\mu^{-1}(U))$ and by~\eqref{property_V}, $$\E[F_{\nu}^{-1}(V)|F^{-1}_\mu(V)]=\E[\E[F_{\nu}^{-1}(V)|U]|F^{-1}_\mu(V)]=\E[\cT(F_\mu^{-1}(V))|F^{-1}_\mu(V)]=\cT(F_\mu^{-1}(V)).$$ Hence the single element of $\tilde{{\mathcal I}}^\nu_\mu$  is the law ${\mathcal T}\#\mu$ of $\cT(F_\mu^{-1}(V))$. Since $\cT$ is nondecreasing, $\cT(F_\mu^{-1}(V))=F_{{\mathcal T}\#\mu}^{-1}(V)$ a.s. and $\E[F_{\nu}^{-1}(V)|F^{-1}_{{\mathcal T}\#\mu}(V)]=F_{{\mathcal T}\#\mu}^{-1}(V)$ a.s.. Hence the law of $(F_{{\mathcal T}\#\mu}^{-1}(V),F_{\nu}^{-1}(V))$, which is the single element of $\Pi^{opt}({\mathcal T}\#\mu,\nu)$, is a martingale coupling. Since all the martingale couplings share the quadratic cost $\int_{\R}y^2\nu(dy)-\int_{\R}({\mathcal T}(x))^2\mu(dx)$, each martingale coupling belongs to $\Pi^{opt}({\mathcal T}\#\mu,\nu)$ and is therefore equal to the previous one.
\end{adem}

In dimension $d=1$, there is a single element $\eta\in\tilde{\mathcal I}^\nu_\mu$, a unique element in $\Pi^{opt}(\mu,\eta)$ and the unique martingale coupling between $\eta$ and $\nu$ is $W_2$-optimal. We now provide an example in dimension~$d=2$ where these properties fail. \begin{ax}
Let $\mu=\frac{1}{2}\left(\delta_{(-1,0)}+\delta_{(1,0)}\right)$ and
  $\nu=\frac{1}{2}\left(\delta_{(0,-1)}+\delta_{(0,1)}\right)$. Since $|(0,-1)-(-1,0)|=|(0,1)-(-1,0)|=|(0,-1)-(1,0)|=|(0,1)-(1,0)|$, any coupling between $\mu$ and $\nu$ is $W_2$-optimal. The couplings write $\mu(dx)k_p(x,dy)$ with $k_p((-1,0),dy)=\left(p\delta_{(0,-1)}+(1-p)\delta_{(0,1)}\right)(dy)$ and $k_p((1,0),dy)=\left((1-p)\delta_{(0,-1)}+p\delta_{(0,1)}\right)(dy)$ for $p\in(0,1)$. One has ${\mathcal T}_p((-1,0))=(0,1-2p)$, ${\mathcal T}_p((1,0))=(0,2p-1)$, and $\eta_p=\frac{1}{2}\left(\delta_{(0,1-2p)}+\delta_{(0,2p-1)}\right)$. Any coupling between $\mu$ and $\eta_p$ is $W_2$-optimal and as soon as $p\neq 1/2$, there is an optimal coupling different from $(I_2,{\mathcal T}_p)\#\mu$. Moreover, unless $p\in\{0,1/2,1\}$, the martingale coupling between $\eta_p$ and $\nu$ is not $W_2$-optimal.
\end{ax}

According to the next theorem, we can find elements $\eta$ in $\tilde{\mathcal I}^\nu_\mu$ such that $\Pi^{opt}(\mu,\eta)=\{(I_d,T)\#\mu\}$ for some measurable transport map $T$ by minimizing over ${\mathcal I}^\nu_\mu$ the integral of a strictly convex function.

\begin{thm}\label{thm_etaphi}
Let $\mu,\nu\in{\mathcal P}_2(\R^d)$, $\phi:\R^d \rightarrow \R$ be strictly convex such that 
$\int_{\R^d}\phi(y)\nu(dy)<\infty$ and ${\mathcal I}_{\mu,\phi}^\nu:=\{\eta\in {\mathcal I}_\mu^\nu:\int_{\R^d} \phi(z) \eta(dz)=\inf_{\eta\in{\mathcal I}_\mu^\nu}\int_{\R^d}\phi(z) \eta(dz)\}$. We have $\emptyset\neq{\mathcal I}_{\mu,\phi}^\nu\subset\tilde{{\mathcal I}}_\mu^\nu$ and for each $\eta\in{\mathcal I}^\nu_{\mu,\phi}$, $\Pi^{opt}(\mu,\eta)=\{(I_d,T)\#\mu\}$ for some measurable transport map $T:\R^d\to\R^d$. Moreover, there is a single $\eta_\phi\in{\mathcal I}_{\mu,\phi}^\nu$ such that $\int_{\R^d}|z|^2\eta_\phi(dz)=\inf_{\eta\in{\mathcal I}_{\mu,\phi}^\nu}\int_{\R^d}|z|^2\eta(dz)$. Last, there is a single element $\underline\eta$ in ${\mathcal I}_{\mu,|x|^2}^\nu$.
\end{thm}
This theorem permits to select extreme elements of ${\mathcal I}_\mu^\nu$ and provides the following characterization of the existence of a minimal element for the convex order in this set. 
\begin{acor}\label{corimin} For $\mu,\nu\in{\mathcal P}_2(\R^d)$, there exists  $\eta_0 \in{\mathcal P}_2(\R^d)$ such that ${\mathcal I}_\mu^\nu=\{ \eta_0 \lecx \eta \lecx \nu \}$ if and only if $\left\{\eta_\phi:\phi:\R^d\to\R^d\mbox{ strictly convex and such that }\int_{\R^d}\phi(y)\nu(dy)<\infty\right\}=\{\underline\eta\}$ and then $\eta_0=\underline\eta$.
\end{acor}
Let us show the corollary before proving the theorem.
\begin{adem}[of Corollary~\ref{corimin}]
  The necessary condition is obvious. Let us show that it is sufficient. It is enough to check that for any $\phi:\R^d \rightarrow \R$  convex such that $\exists C<\infty,\;\forall x\in\R^d,\;|\phi(x)|\le C(1+|x|)$, we have $\forall \eta \in {\mathcal I}_\mu^\nu, \int_{\R^d} \phi(x) \underline{\eta}(dx) \le \int_{\R^d} \phi(x) \eta(dx)$ (see e.g. Lemma A.1~\cite{ACJ-IHP}). For such a function $\phi$ and for $\varepsilon>0$, $\phi_{\varepsilon}(x):= \phi(x)+\varepsilon |x|^2$ is strictly convex and, since $\eta_{\phi_\varepsilon}=\underline\eta$, we have
  $$\forall \eta \in {\mathcal I}_\mu^\nu, \ \int_{\R^d} \phi_{\varepsilon}(x) \underline{\eta}(dx) \le \int_{\R^d} \phi_{\varepsilon}(x) \eta(dx).$$
  We conclude by letting $\varepsilon \rightarrow 0$ using the dominated convergence theorem. \end{adem}

To prove Theorem~\ref{thm_etaphi}, we will need the following Lemma
\begin{alem}\label{lemunifint}
   Let $\nu$ be a probability measure on $\R^d$ such that $\int_{\R^d}|y|\nu(dy)<\infty$ and $\phi:\R^d\to\R$ a convex function such that $\int_{\R^d}\phi(y)\nu(dy)<\infty$. Then the family of probability measures $\{\phi\#\eta:\eta\lecx\nu\}$ is uniformly integrable.
 \end{alem}
 \begin{adem}[of Lemma~\ref{lemunifint}]
   Let us first suppose that $\phi$ is nonnegative. Let $M\in(0,+\infty)$,
   $\eta\lecx\nu$ and $m$ be a martingale kernel such that $\int_{x\in\R^d}\eta(dx)m(x,dy)=\nu(dy)$. Using Jensen's inequality for the first inequality and the Markov inequality combined with $\eta\lecx\nu$ for the third one, we obtain that
   \begin{align*}
     \int_{\R^d}\phi(x)1_{\{\phi(x)\ge M\}}\eta(dx)&\le \int_{\R^d}\int_{\R^d}\phi(y)m(x,dy)1_{\{\phi(x)\ge M\}}\eta(dx)\\&\le \int_{\R^d\times \R^d}\left(\phi(y)1_{\{\phi(y)\ge \sqrt{M}\}}+\sqrt{M}1_{\{\phi(x)\ge M\}}\right)m(x,dy)\eta(dx)\\
                                                   &=\int_{\R^d}\phi(y)1_{\{\phi(y)\ge \sqrt{M}\}}\nu(dy)+\sqrt{M}\int_{\R^d}1_{\{\phi(x)\ge M\}}\eta(dx)\\
                                                     &\le \int_{\R^d}\phi(y)1_{\{\phi(y)\ge \sqrt{M}\}}\nu(dy)+\frac{1}{\sqrt{M}}\int_{\R^d}\phi(y)\nu(dy).
   \end{align*}
   Hence $\lim_{M\to\infty}\sup_{\eta\lecx\nu}\int_{\R^d}\phi(x)1_{\{\phi(x)\ge M\}}\eta(dx)=0$. In particular, the family $\{|x|\#\eta:\eta\lecx\nu\}$ is uniformly integrable. When the sign of $\phi$ is not constant, we obtain a nonnegative convex function $\tilde\phi$ such that $\int_{\R^d}\tilde\phi(y)\nu(dy)<\infty$ by addition to $\phi$ of a suitable affine function $\psi$. The conclusion follows from the uniform integrability of both the families  $\{\psi\#\eta:\eta\lecx\nu\}$ and $\{\tilde\phi\#\eta:\eta\lecx\nu\}$. 
 \end{adem}
\begin{adem}[of Theorem~\ref{thm_etaphi}]
Let $(\eta_n)_{n\in\N}$ be a sequence in ${\mathcal I}_\mu^\nu$ minimizing $\int_{\R^d}\phi(z)\eta(dz)$. For $n\in\N$, let $\mu(dx)q_n(x,dz)\in\Pi^{opt}(\mu,\eta_n)$ and $\eta_n(dz)m_n(z,dy)$ be a martingale coupling between $\eta_n$ and $\nu$. By the second part in Proposition~\ref{lemstructoplan}, $\mu(dx)q_nm_n(x,dy)\in\Pi^{opt}(\mu,\nu)$. Up to extracting a subsequence, we may suppose that $(\mu(dx)q_n(x,dz)m_n(z,dy))_n$  converges weakly to $\mu(dx)r_\infty(x,dz,dy)$ where $\mu(dx)\int_{z\in\R^d}r_\infty(x,dz,dy)\in\Pi^{opt}(\mu,\nu)$. Let ${\mathcal T}_\infty(x)=\int_{\R^d\times\R^d}y r_\infty(x,dz,dy)$ and $\eta_\infty={\mathcal T}_\infty\#\mu$. By the first part of Proposition~\ref{lemstructoplan}, $\eta_\infty\in\tilde{\mathcal I}_\mu^\nu$.
Moreover, by the above weak convergence and the uniform integrability deduced from Lemma~\ref{lemunifint},$$\int_{\R^d\times\R^d\times\R^d}\phi(z)\mu(dx)r_\infty(x,dz,dy)=\lim_{n\to\infty}\int_{\R^d}\phi(z)\eta_n(dz).$$ Taking the limit $n\to\infty$ in the equality
$\int_{\R^d\times\R^d\times\R^d} \varphi(x,z)(y-z)\mu(dx)q_n(x,dz)m_n(z,dy)=0$,
we obtain that $\int_{\R^d\times\R^d\times\R^d} \varphi(x,z)(y-z)\mu(dx)r_\infty(x,dz,dy)=0$ for any continuous and bounded function $\varphi:\R^d\times\R^d\to\R$. Hence, for $(X,Z,Y)$ distributed according to $\mu(dx)r_\infty(x,dz,dy)$, $Z=\E[Y|(X,Z)]$ and ${\mathcal T}_\infty(X)=\E[Y|X]=\E[\E[Y|(X,Z)]|X]=\E[Z|X]$. By using Jensen inequality for the conditional expectation, we get $$\int_{\R^d} \phi(z) \eta_\infty(dz)\le\int_{\R^d\times\R^d\times\R^d}\phi(z)\mu(dx)r_\infty(x,dz,dy)=\lim_{n\to\infty}\int_{\R^d}\phi(z)\eta_n(dz).$$
Thus, $\eta_\infty$ satisfies $\int_{\R^d}\phi(z) \eta_\infty(dz)=\inf_{\eta\in{\mathcal I}_\mu^\nu}\int_{\R^d}\phi(z)\eta(dz)$. Hence ${\mathcal I}_{\mu,\phi}^\nu\neq\emptyset$.

Let $\eta\in{\mathcal I}_{\mu,\phi}^\nu$. We now check that $\eta\in\tilde{{\mathcal I}}_{\mu}^\nu$ and $\Pi^{opt}(\mu,\eta)$ is a singleton. Let $\mu(dx)q(x,dz)\in\Pi^{opt}(\mu,\eta)$ and $\eta(dz)m(z,dy)$ be a martingale coupling between $\eta$ and $\nu$. By the second assertion in Proposition~\ref{lemstructoplan}, $\mu(dx)qm(x,dy)\in\Pi^{opt}(\mu,\nu)$ and, by the first assertion, for ${\mathcal T}(x)=\int_{\R^d}y   qm(x,dy)$, ${\mathcal T}\#\mu\in\tilde{{\mathcal I}}_{\mu}^\nu$. By the martingale property of $m$, ${\mathcal T}(x)=\int_{\R^d}zq(x,dz)$ so that ${\mathcal T}\#\mu\lecx\eta$. Since ${\mathcal T}\#\mu\in{\mathcal I}_{\mu}^\nu$ and $\eta\in{\mathcal I}_{\mu,\phi}^\nu$ implies that $\int_{\R^d} \phi(z) {\mathcal T}\#\mu(dz)\ge\int_{\R^d}\phi(z)\eta(dz)$, we deduce with the strict convexity of $\phi$ that $\eta={\mathcal T}\#\mu$ and $\mu(dx) q(x,dz)=\mu(dx)\delta_{{\mathcal T}(x)}(dz)$. Hence any coupling in $\Pi^{opt}(\mu,\eta)$ is given by a map. By the second statement in Lemma~\ref{cor_utile}, we conclude that this set is a singleton.

By repeating the first argument with $(\phi,{\mathcal I}_\mu^\nu)$ replaced by $(|x|^2,{\mathcal I}_{\mu,\phi}^\nu)$ , we obtain the existence of $\eta_\phi\in{\mathcal I}_{\mu}^\nu$ such that $\int_{\R^d}|z|^2\eta_\phi(dz)\le \inf_{\eta\in{\mathcal I}_{\mu,\phi}^\nu}\int_{\R^d}|z|^2\eta(dz)$. Since the construction also reduces the integral of $\phi$, $\eta_\phi\in{\mathcal I}_{\mu,\phi}^\nu$.

Let us now check that if $\tilde\eta\in{\mathcal I}_{\mu,\phi}^\nu$ is such that $\int_{\R^d}|z|^2\tilde\eta(dz)=\inf_{\eta\in{\mathcal I}_{\mu,\phi}^\nu}\int_{\R^d}|z|^2\eta(dz)$, then $\tilde\eta=\eta_\phi$. By the first statement, $\Pi^{opt}(\mu,\eta_\phi)=\{(I_d,T_\phi)\#\mu\}$ and $\Pi^{opt}(\mu,\tilde\eta)=\{(I_d,\tilde T)\#\mu\}$ for measurable transport maps $T_\phi$ and $\tilde T:\R^d\to\R^d$. One has $\int_{\R^d}|z|^2\eta_\infty(dz)=\int_{\R^d}|z|^2\tilde \eta(dz)$ and therefore, since $\eta_\phi,\tilde\eta\in{\mathcal I}_\mu^\nu$,  $W_2^2(\mu,\eta_\phi)=W_2^2(\mu,\tilde\eta)$.
Let now $\bar \eta=\frac{\eta_\phi+\tilde\eta}{2}$. One has $\int_{\R^d}|z|^2\bar \eta(dz)=\int_{\R^d}|z|^2\eta_\phi(dz)=\int_{\R^d}|z|^2\tilde\eta(dz)$. The coupling $\mu(dx)\frac{1}{2}\left(\delta_{T_\phi(x)}(dz)+\delta_{\tilde T(x)}(dz)\right)$ between $\mu$ and $\bar \eta$ implies that $W_2^2(\mu,\bar \eta)\le W_2^2(\mu,\eta_\phi)=W_2^2(\mu,\tilde\eta)$. Since $\eta_\phi\in{\mathcal I}_\mu^\nu$, we deduce that
$$W_2^2(\mu,\nu)\ge W_2^2(\mu,\bar \eta)+\int_{\R^d}|y|^2\nu(dy)-\int_{\R^d}|z|^2\bar \eta(dz).$$
Moreover,  $\bar \eta\lecx\nu$ and combining a  coupling in $\Pi^{opt}(\mu,\bar\eta)$ with a martingale coupling between $\bar \eta$ and $\nu$, we deduce that the previous inequality is an equality so that $\bar \eta\in{\mathcal I}_\mu^\nu$ and $\mu(dx)\frac{1}{2}\left(\delta_{T_\phi(x)}(dz)+\delta_{\tilde T(x)}(dz)\right)\in\Pi^{opt}(\mu,\bar\eta)$. As $\eta_\phi,\tilde\eta\in{\mathcal I}_{\mu,\phi}^\nu$,  $\int_{\R^d}\phi(z)\bar\eta(dz)=\inf_{\eta\in {\mathcal I}_\mu^\nu}\int_{\R^d}\phi(z)\eta(dz)$ and $\bar\eta\in{\mathcal I}_{\mu,\phi}^\nu$. By the first assertion, $\Pi^{opt}(\mu,\bar\eta)=\{(I_d,\bar T)\#\mu\}$ for some measurable transport map $T:\R^d\to\R$. Therefore $\mu(dx)$ a.e., ${T}_\phi(x)=\tilde{T}(x)$  and $\eta_\phi=\tilde\eta$. For the choice $\phi(x)=|x|^2$, we deduce that ${\mathcal I}_{\mu,|x|^2}^\nu$ is a singleton.
\end{adem}

From the equality $W_2^2(\mu,\nu)=W_2^2(\mu,\eta)+\int_{\R^d}|y|^2\nu(dy)-\int_{\R^d}|z|^2\eta(dz)$ valid for $\eta\in{\mathcal I}^\nu_\mu$, we see that minimizing $\int_{\R^d}|z|^2\eta(dz)$ over ${\mathcal I}^\nu_\mu$ is equivalent to minimizing $W_2^2(\mu,\eta)$. Therefore the probability measure $\underline \eta$ can be seen as the $W_2$-projection of $\mu$ on the set ${\mathcal I}_\mu^\nu$. It is in general different from the $W_2$-projection $\mu_{\underline{\mathcal P}(\nu)}$ of $\mu$ on the set  $\underline{\mathcal P}(\nu):=\{\eta: \eta \lecx \nu \}$, which has been studied recently in dimension $d=1$ by Gozlan et al.~\cite{GRSST} and in general dimension $d$ by Alfonsi et al.~\cite{ACJ-IHP} (who also give an explicit formula for the antiderivative of the quantile function of this projection when $d=1$), Alibert et al.~\cite{ABC}, Gozlan and Juillet~\cite{GoJu} and Backhoff-Veraguas et al.~\cite{B-VBP}. Notice that since ${\mathcal I}_\mu^\nu\subset \underline{\mathcal P}(\nu)$, one always has $W_2(\mu,\mu_{\underline{\mathcal P}(\nu)})\le W_2(\mu,\underline \eta)$. \begin{ax}
   For $\mu$ and $\nu$ the respective uniform distributions on $[0,1]$ and $[0,2]$, we have   ${\mathcal I}_\mu^\nu= \{\nu \}$ and thus $\underline\eta=\nu$. By using the characterization in Theorem 2.6~\cite{ACJ-IHP}, we obtain that the $W_2$-projection $\mu_{\underline{\mathcal P}(\nu)}$ of $\mu$ on the set  $\underline{\mathcal P}(\nu)$ is the uniform distribution on $[1/2,3/2]$.
\end{ax}

The next example shows that the set $$\left\{\eta_\phi:\phi:\R^d\to\R^d\mbox{ strictly convex and such that }\int_{\R^d}\phi(y)\nu(dy)<\infty\right\}$$ may contain distinct elements.
\begin{ax}
   Let $\mu=\frac{1}{2}(\delta_{(-1,0)}+\delta_{(1,0)})$ and $\nu=\frac{1}{4}(\delta_{(-1,-1)}+\delta_{(0,-1)}+\delta_{(0,1)}+\delta_{(1,1)})$. Any optimal coupling between $\mu$ and $\nu$ can be written as $\mu(dx)k_p(x,dy)$ with $k_p((-1,0),dy)=\frac{1}{2}(\delta_{(-1,-1)}+p\delta_{(0,-1)}+(1-p)\delta_{(0,1)})(dy)$ and $k_p((1,0),dy)=\frac{1}{2}(\delta_{(1,1)}+(1-p)\delta_{(0,-1)}+p\delta_{(0,1)})(dy)$ for $p\in[0,1]$. One has ${\mathcal T}_p((-1,0))=(-1/2,-p)$ and ${\mathcal T}_p((1,0))=(1/2,p)$. The measures $\eta_p=\frac{1}{2}(\delta_{(-1/2,-p)}+\delta_{(1/2,p)})$ are not comparable for the convex order since for $p\neq p'$ there is no martingale coupling between $\eta_p$ and $\eta_{p'}$. Moreover, for each $p\in[0,1]$ the unique optimal transport plan $\delta_{((-1,0),(-1/2,-p))}+\delta_{((1,0),(1/2,p))}$ between $\mu$ and $\eta_p$ is given by a map. For this example, $\underline \eta=\eta_0=\frac{1}{2}\left(\delta_{(-1/2,0)}+\delta_{(1/2,0)}\right)$ and $\eta_p=\eta_{\phi_p}$, with $\phi_p(x)=x_1^2+(x_2-2px_1)^2$. The $W_2$-optimal couplings between $\underline\eta$ and $\nu$ can be written as $\eta_0(dz)k_p(2z,dy)$ for $p\in[0,1]$ and in particular the unique martingale coupling $\eta_0(dz)k_0(2z,dy)$ is optimal.
 \end{ax}
 The last example shows that, unlike in the previous one, the martingale couplings between $\underline\eta$ and $\nu$ are not necessarily $W_2$-optimal (even when $\Pi^{opt}(\mu,\nu)$ is a singleton).
\begin{ax}
  Let $\mu=\frac{1}{2}\left(\delta_{(-1,0)}+\delta_{(1,0)}\right)$, $\nu_a=\frac{1}{4}\left(\delta_{(-1,-1)}+\delta_{(-1,2a+1)}+\delta_{(1,-2a-1)}+\delta_{(1,1)}\right)$ with $a\in\R$. The unique $W_2$-optimal coupling between $\mu$ and $\nu_a$ is $\mu(dx)k_a(x,dy)$  with $k_a((-1,0),dy)=\frac{1}{2}(\delta_{(-1,-1)}+\delta_{(-1,2a+1)})(dy)$ and $k_a((1,0),dy)=\frac{1}{2}(\delta_{(1,-2a-1)}+\delta_{(1,1)})(dy)$ so that $\eta_a=\frac{1}{2}\left(\delta_{(-1,a)}+\delta_{(1,-a)}\right)$. Since
  $|(-1,-1)-(-1,a)|^2-|(1,1)-(-1,a)|^2=(a+1)^2-4-(a-1)^2=4(a-1)$, for $a>1$,
\begin{align*}
   W_2^2(\eta_a,\nu_a)&=\frac{1}{2}\left((a+1)^2+4+(a-1)^2\right)<(a+1)^2=\frac{1}{2}\left(3+(2a+1)^2\right)-(1+a^2)\\&=\int|y|^2\nu_a(dy)-\int|z|^2\eta_a(dz),
\end{align*}
so that the martingale coupling between $\eta_a$ and $\nu_a$ is not $W_2$-optimal.
\end{ax}

\section{Differentiability of the squared quadratic Wasserstein distance}\label{secldif}

We now present the notion of differentiability introduced by Lions~\cite{Lions}. Let $f:\mathcal{P}_2(\R^d) \rightarrow \R$. We consider an atomless probability space $(\Omega,\mathcal{A},\P)$ and denote by $L^2(\Omega,\P;\R^d)$ the set of $\R^d$-valued square integrable random variables on this space. The lift of the function~$f$ on $L^2(\Omega,\P;\R^d)$ is the function $F:L^2(\Omega,\P;\R^d)\rightarrow \R$ such that
$$ \forall X \in L^2(\Omega,\P;\R^d), \ F(X)=f(\mathcal{L}(X)),$$
where $\mathcal{L}(X) \in \mathcal{P}_2(\R^d)$ is the probability distribution of $X$. The atomless property is equivalent to the existence of a random variable $U:\Omega \rightarrow \R$ uniformly distributed on $[0,1]$ (see e.g. Proposition A.27 in~\cite{FoSc}). By the fundamental Theorem of simulation (see e.g. Bouleau and Lépingle~\cite{BoLe}, Theorem A.3.1 p. 38), it ensures the existence on $(\Omega,\mathcal{A},\P)$ of a random variable  distributed according to each probability measure on each Polish space,  and in particular of $X:\Omega \rightarrow \R^d$  distributed according to $\mu$, for each $\mu\in{\mathcal P}_2(\R^d)$.

\begin{aadef} A function $f:\mathcal{P}_2(\R^d)\rightarrow \R$ is $L$-differentiable at $\mu\in \mathcal{P}_2(\R^d)$ if there exists $X\in L^2(\Omega,\P;\R^d)$ such that $X\sim\mu$ and $F$ is Fréchet differentiable at $X$. 
\end{aadef}

Let $f:\mathcal{P}_2(\R^d)\rightarrow \R$ and $F(X)= f(\mathcal{L}(X))$ for $X\in L^2(\Omega,\P;\R^d)$. The Fréchet differentiability of~$F$ at~$X$ amounts to the existence of a bounded linear operator $D^F_{X}: L^2(\Omega,\P;\R^d) \rightarrow \R$ such that $F(X+Y)=F(X)+D^F_{X}(Y)+\|Y\|_2 \varepsilon_{X}(Y)$, where $\varepsilon_{X}(Y)\rightarrow 0$ as $\|Y\|_2 \rightarrow 0$. By the Riesz representation theorem, there is a unique $DF(X)\in  L^2(\Omega,\P;\R^d)$ such that $\forall Y\in L^2(\Omega,\P;\R^d),\ D^F_{X}(Y)=\E[DF(X).Y]$, and we will call later on $DF(X)$ the Fréchet derivative of $F$ at $X$. From Theorem~6.2 in~\cite{Cardaliaguet}, if $f$ is $L$-differentiable at $\mu\in \mathcal{P}_2(\R^d)$, then $F$ is Fréchet differentiable at $X$ for all $X\in L^2(\Omega,\P;\R^d)$ such that $\mu=\mathcal{L}(X)$. Besides, the law of $(X,DF(X))$ does not depend on $X$ by Proposition 5.24 in~\cite{CaDe}.
According to Theorem 6.5~\cite{Cardaliaguet}, under the additional continuous differentiability assumption, the Fréchet derivative $DF(X)$ is equal to $g(X)$ for some measurable function $g$. According to Corollary 3.22~\cite{GaTu}, the  continuous differentiability assumption is not needed. We will provide a simple and direct proof of this result, see Lemma~\ref{lemdiffoncx}.

We now state the main result of this section that  characterizes the differentiability of the square quadratic Wasserstein distance. To do so, we first exhibit the lift of  the Wasserstein distance. Let $\mu,\nu \in \mathcal{P}_2(\R^d)$. From the atomless property, there exist random variables $X\sim \mu$ and $Y\sim\nu$ on $(\Omega, \mathcal{A}, \P)$.
The dual formulation (see for instance Theorem 5.10 in~\cite{Villani})
\begin{equation}\label{lifted_wasserstein}
  W_2^2(\mu,\nu)=\sup_{\psi\in L^1(\mu),\tilde \psi\in L^1(\nu):\psi(x)+\tilde \psi(y)\le |x-y|^2}\E\left[\psi(X)+\tilde\psi(Y)\right]=:\mathbb{W}_2^2(X,Y)
\end{equation}
permits to lift $W_2^2$ to $L^2(\Omega,\P;\R^d)$. 
\begin{thm}\label{thm_dif}
   For $\mu,\nu\in{\mathcal P}_2(\R^d)$, the mapping ${\mathcal P}_2(\R^d)\ni\sigma\mapsto W_2^2(\sigma,\nu)$ is $L$-differentiable at $\mu$ iff there exists a measurable function $T:\R^d\to\R^d$ such that $\Pi^{opt}(\mu,\nu)=\{(I_d,T)\#\mu\}$ and then the Fréchet derivative of the function $Z\mapsto \mathbb{W}_2^2(Z,Y)$ at $X\sim \mu$ is given by $2(X-T(X))$.
\end{thm}
\begin{arem}
\begin{itemize}
\item In particular, since the only coupling $\pi \in \Pi^{opt}(\nu,\nu)$  is $(I_d,I_d)\#\nu$, ${\mathcal P}_2(\R^d)\ni\sigma\mapsto W_2^2(\sigma,\nu)$ is differentiable at $\nu$ with a vanishing Fréchet derivative.
  \item According to Proposition \ref{propnondifvois}, if $\nu$ is not a Dirac mass, then there is no $\mu\in{\mathcal P}_2(\R^d)$ such that ${\mathcal P}_2(\R^d)\ni\sigma\mapsto W_2^2(\sigma,\nu)$ is  differentiable in a neighbourhood of $\mu$.
\end{itemize}
\end{arem}

The $L$-differentiability is equivalent to the geometric differentiability (\cite{GaTu}, Corollary 3.22). As explained in the introduction, the sufficient condition in Theorem~\ref{thm_dif}  can be deduced from this equivalence, Theorem 10.2.6~\cite{AGS} and Proposition 4.3~\cite{AmGa}.

We are going to give a probabilistic proof of Theorem~\ref{thm_dif} by working with the $L$-differentiability.  The two following lemmas are needed: the first one is used to get the necessary condition  while the second is used for the sufficient condition. Their proofs are postponed after the proof of the theorem.

\begin{alem}\label{lemdiffoncx}Let $F:L^2(\Omega,\P;\R^d)\to\R$ be law invariant. If $F$ is Fréchet differentiable at $X\sim\mu$, then its Fréchet derivative is equal to $g(X)$ for some measurable function $g\in L^2(\R^d,\mu;\R^d)$ and it is differentiable with Fréchet derivative $g(\tilde X)$ at each $\tilde X\sim\mu$ in $L^2(\Omega,\P;\R^d)$.
\end{alem}
Let us note that this result is also a consequence of the work by Gangbo and Tudorascu~\cite{GaTu}. Here, we provide an alternative simple probabilistic proof of this fact. 
Wu and Zhang (Proposition 1,~\cite{WuZhang}) already gave a different probabilistic proof when $X$ is discrete.

\begin{alem}\label{lemcvloi}
   Let $\mu,\nu\in{\mathcal P}_2(\R^d)$ be such that there exists $T:\R^d \rightarrow \R^d$ measurable such that $\Pi^{opt}(\mu,\nu)=\{(I_d,T)\#\mu\}$. Let also $(\mu_n)_n$ be a sequence of elements of ${\mathcal P}_2(\R^d)$ converging weakly to $\mu$ and such $\lim_{n\to\infty}W_2(\mu_n,\nu)=W_2(\mu,\nu)$. If (on a single probability space), $X\sim\mu$ and for $n\in\N$, $(X_n,Y_n)$ is such that $X_n\sim\mu_n$, $Y_n\sim\nu$, $W_2^2(\mu_n,\nu)=\E\left[|X_n-Y_n|^2\right]$ and $X_n\stackrel{\rm Pr}{\longrightarrow} X$ as $n\to\infty$, then $$\lim_{n\to\infty}\E\left[|X_n-X|^2+|Y_n-T(X)|^2\right]=0.$$
\end{alem}
\begin{arem}
   The fact that $\lim_{n\to\infty}\E\left[|X_n-X|^2\right]=0$ implies that $\lim_{n\to\infty}W_2(\mu_n,\mu)=0$. \end{arem}
\begin{adem}[of Theorem~\ref{thm_dif}]
Let us first assume $\Pi^{opt}(\mu,\nu)\neq\{(I_d,T)\#\mu\}$. The existence on the lifted probability space of a random variable uniformly distributed on $[0,1]$ combined with~\cite{BoLe} Theorem A.3.1. and Lemma~\ref{cor_utile} ensures the existence on this space of $(X,Y)$ with $X\sim\mu$, $Y\sim\nu$, $W_2^2(\mu,\nu)=\E[|Y-X|^2]$ and $\E[|Y-\E[Y|X]|^2]>0$. Let $\xi=Y-\E[Y|X]$. One has $\E[\xi|X]=0$ a.s. so that for $h:\R^d\to\R^d$ measurable and such that $h(X)$ is square integrable, \begin{equation}
   \E[h(X).\xi]=\E[h(X).\E[\xi|X]]=0.\label{moynuld}
\end{equation} On the other hand, denoting by $\mu_n$ the distribution of $X+\xi_n$ where $\xi_n=\frac{\xi}{n}$, we have
\begin{align*}
   \mathbb{W}_2^2(X+\xi_n,Y)=W_2^2(\mu_n,\nu)\le \E[|X+\xi_n-Y|^2]&=\E[|X-Y|^2]+\frac{2}{n}\E[(X-Y).\xi]+\frac{\E[|\xi|^2]}{n^2}\\
&=W_2^2(\mu,\nu)-\frac{2}{n}\E[Y.\xi]+\frac{\E[|\xi|^2]}{n^2}\\&=\mathbb{W}_2^2(X,Y)-\left(\frac{2}{n}-\frac{1}{n^2}\right)\E[|Y-\E[Y|X]|^2],
\end{align*}
where we used \eqref{moynuld} for the second equality and the definition of $\xi$ for the third. If $\sigma\mapsto W_2^2(\sigma,\nu)$ were $L$-differentiable at $\mu$, then \eqref{moynuld} combined with Lemma~\ref{lemdiffoncx} would imply that, as $n\to\infty$, $$\mathbb{W}_2^2(X+\xi_n,Y)-\mathbb{W}_2^2(X,Y)=o(\|\xi_n\|_2),$$  which does not hold since $\|\xi_n\|_2=\frac{\E^{1/2}[|Y-\E[Y|X]|^2]}{n}$.
  
   Now, we assume that $\Pi^{opt}(\mu,\nu)=\{(I_d,T)\#\mu\}$ for some measurable transport map $T:\R^d\to\R^d$.  Let, on the lifted probability space, $X\sim\mu$, $Y\sim \nu$ and $(\xi_n)_n$ be a sequence of square integrable $\R^d$-valued random vectors such that $\|\xi_n\|_2:=\E^{1/2}\left[|\xi_n|^2\right]$ tends to $0$ as $n\to\infty$. We denote by $\mu_n$ the law of $X+\xi_n$. Let $Y_n\sim\nu$ such that $W_2^2(\mu_n,\nu)=\E[|X+\xi_n-Y_n|^2]$ be defined on a possible enlargement of the lifted probability space. We have
\begin{align*}
   W_2^2(\mu_n,\nu)&\le \E[|X+\xi_n-T(X)|^2]=\E[|X-T(X)|^2]+2\E[(X-T(X)).\xi_n]+\E[|\xi_n|^2]\\&=W_2^2(\mu,\nu)+2\E[(X-T(X)).\xi_n]+\E[|\xi_n|^2].
\end{align*}
On the other hand,
\begin{align*}
   W_2^2(\mu,\nu)&\le \E[|X-Y_n|^2]=\E[|X+\xi_n-Y_n|^2]-2\E[(X-Y_n).\xi_n]-\E[|\xi_n|^2]\\&=W_2^2(\mu_n,\nu)-2\E[(X-T(X)).\xi_n]-\E[|\xi_n|^2]+2\E[(Y_n-T(X)).\xi_n].
\end{align*}
With Cauchy-Schwarz inequality, we deduce that
\begin{align*}
   |W_2^2(\mu_n,\nu)-W_2^2(\mu,\nu)-2\E[(X-T(X)).\xi_n]|\le \|\xi_n\|_2\left(\|\xi_n\|_2+\|Y_n-T(X)\|_2\right).
\end{align*}
Note that from~\eqref{lifted_wasserstein} the left-hand side is equal to $|\mathbb{W}_2^2(X+\xi_n,Y)-\mathbb{W}_2^2(X,Y)-2\E[(X-T(X)).\xi_n]|$ and is thus well defined on the original lifted probability space, as required by the definition of the Lions derivative. Now, Lemma~\ref{lemcvloi} applied with $X_n=X+\xi_n$ ensures that $\lim_{n\to \infty}\|Y_n-T(X)\|_2 = 0$ so that $\sigma\mapsto W_2^2(\sigma,\nu)$ is $L$-differentiable at $\mu$ with $\partial_\mu W_2^2(\mu,\nu)(x)=2(x-T(x))$.   
\end{adem}

\begin{adem}[of Lemma~\ref{lemdiffoncx}]
   By the fundamental theorem of simulation (see e.g. Bouleau and Lépingle~\cite{BoLe}, Theorem A.3.1 p. 38), since the lifted probability space supports a random variable with uniform distribution on $[0,1]$, it also supports a couple $(\tilde X,U)$ with $\tilde X\sim\mu$ and $U$ an independent random variable uniformly distributed on $[0,1]$.  Since $F$ is Fréchet differentiable at $X\sim\mu$,  by Theorem 6.2~\cite{Cardaliaguet} it is also Fréchet differentiable at $\tilde X$.  Let for $i\in\{1,\hdots,d\}$, $DF(\tilde X)_i$ denote the $i$-th coordinate of $DF(\tilde X)$ and $P_i(x,dz,du)$ with respective marginals $Q_i(x,dz)$ and $R_i(x,du)$ denote a regular version of the conditional law of $(DF(\tilde X)_i,U)$ given $\tilde X=x$. Let $g_i(x)=\inf\{z\in\R:Q_i(x,(-\infty,z])\ge 1/2\}$ be the median of $Q_i(x,dz)$. Notice that, by property of the median, $\E[|DF(\tilde X)_i-g_i(\tilde X)||\tilde X]\le\E[|DF(\tilde X)_i-\E[DF(\tilde X)_i|\tilde X]||\tilde X]$ so that
\begin{align*}
   \E[|g_i(\tilde X)|]&\le \E[|DF(\tilde X)_i|]+\E[|DF(\tilde X)_i-g_i(\tilde X)|]\\&\le\E[|DF(\tilde X)_i|]+\E[|DF(\tilde X)_i-\E[DF(\tilde X)_i|\tilde X]|]<\infty.
\end{align*}

Let $$v^{\pm}_i(x)=\inf\{u\in[0,1]:P_i(x,\{g_i(x)\}\times [0,u])\ge (Q_i(x,(-\infty,g_i(x)))-Q_i(x,(g_i(x),+\infty)))^\pm\}.$$ By independence of $\tilde X$ and $U$, there is a Borel subset $A$ of $\R^d$ with $\mu(A)=0$ such that for $x\notin A$, $R_i(x,du)$ is the Lebesgue measure on $[0,1]$. Since $$Q_i(x,(-\infty,g_i(x)))\vee Q_i(x,(g_i(x),+\infty))\le \frac{1}{2}\le Q_i(x,(-\infty,g_i(x)))\wedge Q_i(x,(g_i(x),+\infty))+Q_i(x,\{g_i(x)\}),$$ for $x\notin A$, $P_i(x,\{g_i(x)\}\times [0,v^\pm_i(x)])=(Q_i(x,(-\infty,g_i(x)))-Q_i(x,(g_i(x),+\infty)))^\pm$.
   
The random variables $\xi^i_+=1_{\{DF(\tilde X)_i>g_i(\tilde X)\}}+1_{\{DF(\tilde X)_i=g_i(\tilde X),U\le v_i^+(\tilde X)\}}$ and $\xi^i_-=1_{\{DF(\tilde X)_i< g_i(\tilde X)\}}+1_{\{DF(\tilde X)_i=g_i(\tilde X),U\le v_i^-(\tilde X)\}}$ are such that $(\tilde X,\xi^i_+)$ and $(\tilde X,\xi^i_-)$ have the same distribution: indeed, conditionally on $\tilde{X}=x$, these are Bernoulli random variables of parameter $Q_i(x,(-\infty,g_i(x)))\vee Q_i(x,(g_i(x),+\infty))$. Therefore $\E[g_i(\tilde X)\xi^i_+]=\E[g_i(\tilde X)\xi^i_-]$ and, denoting by $e_i$ the $i$-th vector of the canonical basis of $\R^d$, for each $\varepsilon\in[0,1]$, $\tilde X+\varepsilon\xi^i_+e_i$ and $\tilde X+\varepsilon\xi^i_-e_i$ have the same distribution so that $F(\tilde X+\varepsilon\xi^i_+e_i)=F(\tilde X+\varepsilon\xi^i_-e_i)$. Hence $\E(\xi^i_+DF(\tilde X)_i)=\E(\xi_-^iDF(\tilde X)_i)$. We deduce that $$0=\E[(DF(\tilde X)_i-g_i(\tilde X))(\xi^i_+-\xi^i_-)]=\E[|DF(\tilde X)_i-g_i(\tilde X)|]$$
and conclude that $\P\left(DF(\tilde X)=g(\tilde X)\right)=1$.
Proposition 5.24~\cite{CaDe} ensures that the couples $(X,DF(X))$ and $(\tilde X,DF(\tilde X))$ share the same distribution and therefore $\P\left(DF(X)=g(X)\right)=1$.
 \end{adem}

\begin{adem}[of Lemma~\ref{lemcvloi}]
   Let $\eta_n$ and $\eta^n_{34}$ respectively denote the distributions of $(X,T(X),X_n,Y_n)$ and $(X_n,Y_n)$. Since $(\mu_n)_n$ converges weakly to $\mu$, this sequence is tight and we deduce that $(\eta_n)_n$ is tight. Let us consider a subsequence weakly converging to $\eta^\infty$ and that we still index by $n$ for notational simplicity. From the convergence $X_n\stackrel{\rm Pr}{\longrightarrow} X$ as $n\to\infty$, we deduce that $(X,T(X),X_n)\stackrel{\rm Pr}{\longrightarrow} (X,T(X),X)$. Hence the marginal $\eta^\infty_{123}$ of the triplet of the three first coordinates under $\eta^\infty$ is $\eta^\infty_{123}=(I_d,T,I_d)\#\mu$. Next, the marginal $\eta^\infty_{34}$ of the couple of the two last coordinates is a coupling between $\mu$ and $\nu$ such that $\int_{\R^d\times\R^d}|x-y|^2\eta^\infty_{34}(dx,dy)\le\liminf_{n\to\infty}\int_{\R^d\times\R^d}|x-y|^2\eta^n_{34}(dx,dy)$. Since $\int_{\R^d\times\R^d}|x-y|^2\eta^n_{34}(dx,dy)=\E\left[|X_n-Y_n|^2\right]=W_2^2(\mu_n,\nu)$ converges to $W_2(\mu,\nu)$ as $n\to\infty$, $\eta^\infty_{34}\in\Pi^{opt}(\mu,\nu)$. Therefore $\eta^\infty_{34}=(I_d,T)\#\mu$ and $\mu(dw)$ a.e. the conditional law  of the fourth coordinate given that the third is equal to $w$ is $\delta_{T(w)}(dz)$. Since under $\eta^\infty_{123}$ the two first coordinates are a function of the third one, this is also the conditional law of the fourth coordinate given that the three first coordinates are equal to $(w,T(w),w)$. Hence 
$\eta^\infty(dx,dy,dw,dz)=\mu(dw)\delta_{(w,T(w))}(dx,dy)\delta_{T(w)}(dz)$ so that $\eta^\infty=(I_d,T,I_d,T)\#\mu$. Since the weak limit does not depend on the subsequence, the whole sequence $(\eta_n)_n$ converges weakly to $(I_d,T,I_d,T)\#\mu$.\\
Since $|Y_n-T(X)|\le 2|Y_n|1_{\{|Y_n|\ge|T(X)|\}}+2|T(X)|1_{\{|Y_n|<|T(X)|\}}$, for $m>0$, we have $$|Y_n-T(X)|^21_{\{|Y_n-T(X)|\ge m\}}\le 4|Y_n|^21_{\{|Y_n|\ge m/2\}}+4|T(X)|^21_{\{|T(X)\ge m/2\}}.$$ Hence 
$$\E\left[|Y_n-T(X)|^21_{\{|Y_n-T(X)|\ge m\}}\right]\le 8\int_{\R^d}|x|^21_{\{|x|\ge m/2\}}\nu(dx)$$
which provides the uniform integrability needed to conclude that $\lim_{n\to\infty}\E\left[|Y_n-T(X)|^2\right]=0$. The convergence of $\E\left[|X_n-Y_n|^2\right]=W_2^2(\mu_n,\nu)$ to $W_2^2(\mu,\nu)=\E\left[|X-T(X)|^2\right]$ as $n\to\infty$ together with the convergence in probability of $(X_n,Y_n)$ to $(X,T(X))$ implies that $$\lim_{n\to\infty}\E\left[||X-T(X)|^2-|X_n-Y_n|^2|\right]=0$$ and therefore that the random variables $(|X_n-Y_n|^2)_n$ are uniformly integrable. From the inequality $|X_n-X|^2\le 3(|X_n-Y_n|^2+|Y_n|^2+|X|^2)$ and the convergence of $(Y_n)_n$ to $T(X)$ in quadratic mean, we deduce that the random variables $(|X_n-X|^2)_n$ are uniformly integrable. With the convergence in probability of $(X_n)_n$ to $X$, we conclude that $\lim_{n\to\infty}\E\left[|X_n-X|^2\right]=0$.
\end{adem}

Lemma~\ref{lemdiffoncx} is also useful to check that the notion of $L$-differentiability does not depend on the choice of the atomless lifted probability space~$(\Omega,\mathcal{A},\P)$. We will say that a function $f:\mathcal{P}_2(\R^d)\rightarrow \R$ is {\it $L_\Omega$-differentiable} at $\mu\in \mathcal{P}_2(\R^d)$ if there exists $X\in L^2(\Omega,\P;\R^d)$ such that $\mu=\mathcal{L}(X)$ and $F_\Omega(X)=f(\mathcal{L}(X))$ is Fréchet differentiable at~$X$. The function $F_\Omega$ is called the lift of~$f$ on the probability space $(\Omega,\mathcal{A},\P)$.

\begin{aprop}\label{prop_cle} Let $(\Omega,\mathcal{A},\P)$ and $(\tilde \Omega,\tilde{\mathcal{A}},\tilde{\P})$ be two atomless probability spaces. The function $f:\mathcal{P}_2(\R^d)\rightarrow \R$ is $L_\Omega$-differentiable at $\mu\in \mathcal{P}_2(\R^d)$ iff it is  $L_{\tilde{\Omega}}$-differentiable at $\mu$. 
\end{aprop}

\begin{adem}
  By symmetry, it is enough to prove that the $L_\Omega$-differentiability implies the $L_{\tilde{\Omega}}$-differentiability and the Fréchet derivatives are given by the same function in $L^2(\R^d,\mu;\R^d)$. Let us assume that $f$ is $L_\Omega$-differentiable at $\mu$. The atomless property and the fundamental theorem of simulation ensure the existence on the original lifted space  $(\Omega,\mathcal{A},\P)$ of random variables  $(U,  X)$ such that $U$ is uniformly distributed on $[0,1]$ and independent from $  X \sim \mu$. Then, $F_\Omega$ is Fréchet differentiable at $  X$ and there exists a measurable function $g\in L^2(\R^d,\mu;\R^d)$ such that $DF_\Omega({X})=g({X})$ by Lemma~\ref{lemdiffoncx}. We consider $F_{\tilde{\Omega}}: L^2(\tilde \Omega,\tilde \P;\R^d) \rightarrow \R$ the lift of~$f$ on an atomless  probability space  $(\tilde \Omega,\tilde{\mathcal{A}},\tilde \P)$ and  $\tilde X \sim \mu$ under $\tilde{\P}$.
  Let  $\tilde{Y}\in L^2(\tilde \Omega,\tilde \P;\R^d)$ and $R(x,dy)$ denote a regular version of the conditional law of $\tilde{Y}$ given $\tilde X=x$.  By Lemma~2.22~\cite{Kallenberg}, there exists a measurable function $\rho:\R^d\times[0,1]\to \R^d$ such that for all $x\in \R^d$, $\rho(x,U)$ is  distributed according to $R(x,dy)$. Then, $ Y=\rho(X,U)$ is such that $( X, Y)$ has the same law under~$\P$ as $(\tilde X ,\tilde{Y})$ under $\tilde{\P}$, and therefore $Y\in L^2(\Omega,\P;\R^d)$. We then have  
  \begin{align*}
    F_{\tilde{\Omega}}(\tilde X+ \tilde Y)- F_{\tilde{\Omega}}(X)-\tilde{\E}[ g(\tilde X). \tilde Y]&= F_\Omega( X+  Y)-F_\Omega( X)-\E[g( X).  Y]\\&
    = F_\Omega( X+  Y)-F_\Omega( X)-\E[DF_\Omega( X) . Y]. 
  \end{align*}
With $\E[| Y|^2]=\tilde\E[|\tilde Y|^2]$, we deduce that the Fréchet differentiability of $ F_\Omega$ at $ X$ implies the Fréchet differentiability of $F_{\tilde{\Omega}}$ at $\tilde X$ and $DF_{\tilde{\Omega}}[\tilde X]=g(\tilde X)$.
\end{adem}

By considering the atomless probability space $(B_d,\mathcal{B}(B_d),\textup{Leb})$ where $B_d$ the ball centered at the origin of unit volume in~$\R^d$ endowed with the Borel sets and the Lebesgue measure, Gangbo and Tudorascu~(\cite{GaTu}, Corollary 3.22) have proved the equivalence between the $L_{B_d}$-differentiability (and therefore the $L$-differentiability) and the geometric differentiability, as well as the relation between the two derivatives.

As pointed by one of the referees, one can also deduce Proposition~\ref{prop_cle} from Corollary~3.22~\cite{GaTu} by Gangbo and Tudorascu: $f$ is geometrically differentiable at~$\mu$ if, and only if, $F_{B_d}$ is Fréchet differentiable at any $X\sim \mu$, and in this case we have $DF_{B_d}(X)=\nabla_\mu f(X)$, where $\nabla_\mu f: \R^d \to \R^d$ is the Wasserstein gradient of $f$ at~$\mu$. Let $(\Omega, \mathcal{A}, \P)$ be an atomless probability space and $F_\Omega$ denote the corresponding lift of~$f$. Then, there exists an almost isomorphism $i:B_d\rightarrow \Omega$ that pushes forward the Lebesgue measure to $\P$. We have  $F_{B_d}(X\circ i)=F_\Omega(X)$ for $X\in L^2(\Omega,\P ;\R^d)$. If $X\sim \mu$, then $X\circ i \sim \mu$ and one easily deduces that $F_\Omega$ is Fréchet differentiable at~$X$ with  $DF_\Omega(X)= DF_{B_d}(X\circ i)\circ i^{-1}=\nabla_\mu f(X\circ i )\circ i^{-1} =\nabla_\mu f(X)$. In a symmetric way, if $F_\Omega$ is Fréchet differentiable at $X\sim \mu$, then $F_{B_d}$ is Fréchet differentiable at $X\circ i$ and $DF_{B_d}(X\circ i)= DF_\Omega(X)\circ i=\nabla_\mu f(X\circ i)$, so that $DF_\Omega(X)=\nabla_\mu f(X)$.


\bibliographystyle{plain}
\bibliography{Biblio}

\end{document}